\documentclass[12pt]{article}

\usepackage{amsmath,amssymb,amsthm}
\usepackage {graphicx,authblk}

\newtheorem{theorem}{Theorem}[section]
\newtheorem{lemma}[theorem]{Lemma}
\newtheorem{remark}[theorem]{Remark}
\newtheorem{corollary}[theorem]{Corollary}

\newtheorem {definition}[theorem]{Definition}
\newtheorem{problem}[theorem]{Problem}
\def \RR{\mathbb R}
\def \NN{\mathbb N} 
\def\({\left(}  \def\){\right)}
\def\[{\left[}  \def\]{\right]}

\def \beq {\begin {equation}}
\def \eeq {\end{equation}}
\def \OL {\overline}

\begin {document}
\begin {center}
\Large{Min-Max Polarization for Certain Classes\\
of Sharp Configurations on the Sphere}

\bigskip

Sergiy Borodachov
\date{

\today
}

\bigskip

\normalsize{\it Department of Mathematics, Towson University, 8000 York Rd, Towson, MD, 21252}
\footnote {phone: 410-704-2072; fax 410-704-4149 \\ e-mail: sborodachov@towson.edu}
\end {center}
\large {

\begin {abstract}
We consider the problem of finding an $N$-point configuration on the sphere $S^d\subset \RR^{d+1}$ with the smallest absolute maximum value over $S^d$ of its total potential. The potential induced by each point ${\bf y}$ in a given configuration at a point ${\bf x}\in S^d$ is $f\(\left|{\bf x}-{\bf y}\right|^2\)$, where $f$ is continuous on $[0,4]$ and completely monotone on $(0,4]$, and $\left|{\bf x}-{\bf y}\right|$ is the Euclidean distance between points~${\bf x}$ and ${\bf y}$. We show that any sharp point configuration $\OL\omega_N$ on $S^d$, which is antipodal or is a spherical design of an even strength is a solution to this problem. We also prove that the absolute maximum over $S^d$ of the potential of any such configuration $\OL\omega_N$ is attained at points of $\OL\omega_N$.
\end {abstract}

\normalsize{

\smallskip

}

\section {Introduction and setting of the problem}

Let $S^{d}:=\{(x_1,\ldots,x_{d+1})\in \RR^{d+1} : x_1^2+\ldots+x_{d+1}^2=1\}$ denote the unit sphere in the $(d+1)$-dimensional Euclidean space. The classical discrete polarization problem requires finding positions of $N$ points on $S^d$ with the largest absolute minimum over $S^{d}$ of their total potential. It was studied in works \cite {Sto1975circle,Amb2009,AmbBalErd2013,ErdSaf2013,HarKenSaf2013,BorBos2014,Su2014,BosRue2017,BorHarRez2018,AmbNie2019,Borsimplex}. In particular, for the Riesz and logarithmic potential on the sphere, the complete solution to the discrete polarization problem is known only for $N$ points on the unit circle $S^1$ and for up to $d+2$ points on~$S^d$, $d\geq 2$. The continuous version of the polarization problem has also been studied (for references, cf., e.g., \cite {BorHarSafbook}). 

For potential functions bounded above, some authors also dealt with the problem of minimizing the maximum value of the potential, see \cite {Sto1975circle,BosRue2017,AmbNie2019}. We show here the optimality for this problem of a certain class of sharp spherical configurations.

The problem of finding the absolute minimum or maximum of the total potential of regular point configurations on $S^{d}$ is useful for solving a number of important cases of these polarization problems. Several authors have found absolute extrema for some of these configurations, see \cite {Sto1975circle,Sto1975,NikRaf2011,NikRaf2013}. 
We propose here a general approach to characterizing absolute maxima of the potential of tight spherical designs using the classical Delsarte-Yudin method. For a detailed description of this method, see \cite {CohKum2007,BorHarSafbook} and references therein. 

Let $f:[0,4]\to\RR$ be a function continuous on the interval $[0,4]$. We will call $f$ {\it the potential function} and will specify additional assumptions on $f$ later. The interaction between points ${\bf x}$ and ${\bf y}$ on $S^{d}$ is described by the kernel of the form $K({\bf x},{\bf y})=f(\left|{\bf x}-{\bf y}\right|^2)$. For a given configuration $\omega_N=\{{\bf x}_1,\ldots,{\bf x}_N\}\subset S^d$, denote {\it the potential of $\omega_N$} by
\begin {equation}\label {p_s}
p_f({\bf x},\omega_N):=\sum\limits_{i=1}^{N}f(\left|{\bf x}-{\bf x}_i\right|^2),\ \ \ {\bf x}\in S^d.
\end {equation}

\begin {problem}\label {P2}
{\rm
Find the quantity
\begin {equation}\label {max1'}
Q_f(\omega_N;S^d):=\max\limits_{{\bf x}\in S^d}p_f({\bf x},\omega_N),
\end {equation}
and determine all the points ${\bf x}^\ast\in S^d$, for which the maximum is attained in~\eqref {max1'}.
}
\end {problem}
The classical polarization problem requires (for a given $N$) maximizing the quantity
\begin {equation}\label {P}
P_f(\omega_N;S^{d}):=\min\limits_{{\bf x}\in S^{d}}p_f({\bf x},\omega_N)
\end {equation}
over all $N$-point configurations\footnote{Points in $\omega_N$ can coincide; i.e., it is assumed to be a multiset.} $\omega_N\subset S^{d}$. Since we assumed $f$ to be finite and continuous on $[0,4]$, it makes sense to consider the dual problem. 
\begin {problem}\label {mp}
{\rm
Given $N\in \NN$, find the quantity
\begin {equation}\label {dual}
\mathcal Q_f(S^{d},N):=\min\limits_{\omega_N\subset S^{d}} Q_f(\omega_N;S^{d})
\end {equation}
and $N$-point configurations $\omega_N^\ast\subset S^{d}$ that attain the minimum on the right-hand side of \eqref {dual}.}
\end {problem}
A point configuration $\omega_N$ is a solution to Problem \ref {mp} if and only if it is a solution to the classical polarization problem with the potential function $-f$.

Among kernels of the form $K({\bf x},{\bf y})=f(\left|{\bf x}-{\bf y}\right|^2)$, the most interesting to us are the ones with a completely monotone potential function $f$ on $(0,4]$ or with $f$ that becomes completely monotone after adding an appropriate constant. We recall that an infinitely differentiable function $f:I\to \RR$ is {\it completely monotone} on an interval $I$ if $(-1)^k f^{(k)}$ is non-negative on $I$ for all $k\geq 0$ and {\it strictly completely monotone} if $f$ is completely monotone on $I$ with $(-1)^k f^{(k)}$ being strictly positive in the interior of $I$ for all $k\geq 0$. Such are potential functions $f_s(t)=t^{-s/2}$, $s>0$, on $(0,\infty)$ and $f(t)=e^{-at}$, $a>0$, on $(-\infty,\infty)$. They define the Riesz $s$-kernel $K_s({\bf x},{\bf y}):=\left|{\bf x}-{\bf y}\right|^{-s}$ and the Gaussian kernel $K({\bf x},{\bf y})=e^{-a \left|{\bf x}-{\bf y}\right|^2}$, respectively. 
Furthermore, potential functions $f_s(t)=-t^{-s/2}$ for $-2<s<0$ and $f_{\log}(t)=\frac {1}{2}\ln \frac {1}{t}$ become (strictly) completely monotone on $(0,4]$ after adding a certain constant. They give a rise to the Riesz kernel $K_s({\bf x},{\bf y}):=-\left|{\bf x}-{\bf y}\right|^{-s}$ and the logarithmic kernel $K_{\log}({\bf x},{\bf y})=\ln \frac {1}{\left|{\bf x}-{\bf y}\right|}$, respectively.

Any completely monotone $f$ on $(0,4]$ has a {\it convex} derivative of any even order and a {\it concave} derivative of any odd order. It satisfies assumptions of our main results if we assume additionally that it is bounded above and extend it to $t=0$ by the limit value.

Exact solutions to the maximal polarization problem on the sphere are known on $S^1$ for any cardinality $N\geq 3$ and on $S^{d-1}$, $d\geq 3$, for $3\leq N\leq d+1$ (the cases $N=1$ and $N=2$ are trivial). The optimality of the vertices of an equilateral triangle inscribed in $S^1$ was shown for the Riesz potential functions for all $s\neq 0$ by Stolarsky \cite {Sto1975circle} and by Nikolov and Rafailov \cite {NikRaf2011}. The optimality of the vertices of a regular $N$-gon inscribed in $S^1$ was proved for arbitrary $N$ by Ambrus \cite {Amb2009} and Ambrus, Ball, and Erd\'elyi \cite {AmbBalErd2013} (for $s=2$), by Erd\'elyi and Saff \cite {ErdSaf2013} (for $s=4$), and by Hardin, Kendall, and Saff \cite{HarKenSaf2013} for $s\geq -1$, $s\neq 0$, and for $s=\log$. Paper \cite {HarKenSaf2013}, in fact, established this result for any kernel of the form $f(l ({\bf x},{\bf y}))$, where $l({\bf x},{\bf y})$ is the geodesic distance between points ${\bf x}$ and ${\bf y}$ on $S^1$ and $f$ is a decreasing and convex function on $(0,\pi]$ defined by its (possibly infinite) limit value at $t=0$. 

Given $s\in \{-2,-4,\ldots,2-2N\}$, $N\geq 2$, a configuration $\{z_1,\ldots,z_N\}$ of complex numbers on $S^1$ considered as a subset of the complex plane $\mathbb C$ is optimal for the maximal polarization problem on $S^1$ if and only if 
\begin {equation}\label {s/2}
\sum\limits_{i=1}^{N}z_i^k=0,\ \ \ k=1,\ldots,-s/2,
\end {equation}
see the work by Bosuwan and Ruengrot \cite {BosRue2017}. The same is true for the dual problem of minimizing the maximum value of the potential of an $N$-point configuration on $S^1$. Condition \eqref {s/2} holds for the configuration of $N$ equally spaced points on $S^1$. 
The problem of minimizing the quantity $\max_{{\bf x}\in S^1}\sum_{i=1}^{N}\left|{\bf x}\cdot {\bf u}_i\right|^s$ over all configurations $\left\{{\bf u}_1,\ldots,{\bf u}_N\right\}$ on $S^1$ was solved for $0<s\leq 1$ by Ambrus and Nietert in \cite{AmbNie2019} and for $s=2,4,\ldots,2N-2$ in \cite{BosRue2017} (a short proof for the latter case was also given in \cite{AmbNie2019}). 

On the sphere $S^{d-1}$, $d\geq 3$, the solution to the maximal polarization problem has an elementary proof for $3\leq N\leq d$ and any non-increasing and convex potential function $f$ as well as for $s=-2$ and arbitrary cardinality $N\geq 3$. The solution in both cases is any $N$-point configuration with center of mass at the origin (see, e.g., \cite [Section 14.2]{BorHarSafbook}).
For $N=d+1$, the problem turned out to be rather difficult. The optimality of the set of vertices of a regular $d$-simplex $\omega_{d+1}^\ast$ inscribed in $S^{d-1}$ was proved by Su in \cite {Su2014} for $d=3$, $s>-2$, $s\neq 0$, and potential functions
\begin {equation}\label {Rafailov}
f(t)=\begin {cases}\ \ (t+C)^{-s/2}, & s>0, \cr
-(t+C)^{-s/2}, & s<0,\cr
\end {cases} \ \ \ C\geq 0.
\end {equation}
The author of this paper proved in \cite {Borsimplex} a similar result for $d\geq 3$ and potential functions $f$ that are decreasing and convex on $(0,4]$, have a convex or concave derivative $f'$ on~$(0,4)$, and are defined at $t=0$ by their limit value. Another non-trivial result is minimizing
\begin {equation}\label {abs}
\max_{{\bf x}\in S^{d-1}}\sum_{i=1}^{N}\left({\bf x}\cdot {\bf u}_i\right)^2
\end {equation} 
over all configurations $\omega_N=\{{\bf u}_1,\ldots,{\bf u}_N\}\subset S^{d-1}$, see \cite{AmbNie2019}. That paper proved that $\omega_N$ is optimal if and only if vectors ${\bf u}_1,\ldots,{\bf u}_N$ form an isotropic set or, equivalently, a unit-norm tight frame. The optimality of the equitriangular lattice among lattices in $\RR^2$ of the same density was shown by Montgomery \cite {Mon1988} and B\'etermin, Faulhuber, and Steinerberger \cite {BetFauStehex} for min-max and max-min polarization, respectively (for Gaussian kernels).

Exact solutions to the maximal polarization problem on certain sets other than spheres are also known as well as asymptotic results on general classes of sets and the work on the continuous version of the polarization problem that uses potential-theoretic approach. Detailed reviews are given, for example, in book \cite [Section 14.2]{BorHarSafbook} and in paper \cite {Borsimplex}.

In this paper, we solve Problem \ref {mp} for potential functions $f$ continuous on $[0,4]$ and completely monotone on $(0,4]$ (modulo an additive constant) proving the optimality of any strongly sharp or antipodal sharp configuration on $S^d$ (see Definition \ref {Cohn}). We use the universal optimality result by Cohn and Kumar for the discrete minimal energy problem (see \cite {CohKum2007} and references therein to earlier works by Yudin, Kolushov, and Andreev). A detailed description of the Delsarte-Yudin (or linear programming) approach and many references to works using it can be found, for example, in \cite {CohKum2007} and \cite [Chapter 5]{BorHarSafbook}.

The early work by Stolarsky \cite {Sto1975circle,Sto1975} considered Problem \ref {P2} and the problem of finding quantity \eqref {P} for the Riesz potential function $f_s(t)=-t^{-s/2}$ for a certain range of $s<0$ and the vertices of a regular polygon, regular simplex, regular cross-polytope, and a cube as well as the quantity in \eqref {P} for $s>0$ and the vertices of a regular cross-polytope and a cube. Later, Nikolov and Rafailov \cite {NikRaf2011,NikRaf2013} extended these results to all values of $s\neq 0$ and all four types of configurations mentioned above for potential functions
 \eqref {Rafailov}. In all these results, for all four types of configurations, the situation is determined by the sign of the derivative of order $n+1$, where $n$ is the strengh of the configuration as a spherical design.

Concerning general potentials, paper \cite {Borsimplex} showed that $p_f({\bf x},\omega_{d+1}^\ast)$ is minimized at points of $-\omega_{d+1}^\ast$ and maximized at points of $\omega_{d+1}^\ast$ for any potential function $f$ with a concave derivative $f'$ on $(0,4)$. This extended the results of \cite {Sto1975,NikRaf2011} for the regular simplex.
The results of Stolarsky, Nikolov, and Rafailov mentioned above for the maximum of the potential of the vertices of a regular polygon and a regular cross-polytope are extended here
to general potential functions $f$ in a similar way (see Theorems \ref {2m-1w} and~\ref {sharp}). Theorems~\ref {2m-1w} and \ref {sharp} also extend these results to sharp configurations that are antipodal or are designs of even strength. The problem of minimizing the potential of these configurations and of a certain family of non-sharp configurations is studied in \cite {Borsymmetric,Borstiff} using a similar approach. These results were announced in the talk \cite {Bor2022talk}.


\section {Main results}\label {main}

In this section, we state the solution to Problems \ref {P2} and \ref {mp} for the two families of sharp configurations mentioned above (they are known to be tight designs).

For a given ${\bf z}\in S^d$ and a given point configuration $\omega_N=\{{\bf x}_1,\ldots,{\bf x}_N\}\subset S^{d}$,
we denote by 
$$
D({\bf z},\omega_N):=\{{\bf z}\cdot {\bf x}_i : i=1,\ldots,N\}
$$ 
the set of all dot products formed by the point ${\bf z}$ and points from $\omega_N$.

Let $\sigma_{d}$ denote the $d$-dimensional area measure on $S^d$ normalized to be a probability measure; that is, $\sigma_{d}\(\cdot\)=\frac {1}{\mathcal H_{d}(S^{d})}\mathcal H_d|_{S^{d}}\(\cdot\)$, where $\mathcal H_{d}$ is the $d$-dimensional Hausdorff measure in $\RR^{d+1}$.
Following \cite {DelGoeSei1977}, we call a point configuration $\omega_N$ a {\it spherical $n$-design} if, for every polynomial $p$ on $\RR^{d+1}$ of degree at most~$n$,
$$
\frac {1}{N}\sum\limits_{i=1}^{N}p({\bf x}_i)=\int_{S^{d}}p({\bf x})\ \! d\sigma_{d}({\bf x}).
$$
The number $n$ is called the {\it strength} of the spherical design $\omega_N$.

We first consider the class of sharp configuration which are even order designs. The definition of a sharp configuration was first given in \cite {CohKum2007}.
\begin {definition}\label {Cohn}
{\rm
We call a point configuration $\omega_N\subset S^{d}$ {\it $m$-sharp} if $\omega_N$ is a spherical $(2m-1)$-design and there are $m$ distinct values of the dot product formed by any two distinct vectors from $\omega_N$. We call an $m$-sharp configuration {\it strongly $m$-sharp} if it is a spherical $2m$-design. 
}
\end {definition}
As examples of strongly sharp configurations we mention the set of vertices of a regular $(2m+1)$-gon inscribed in $S^1$ (strongly $m$-sharp), the set of $d+2$ vertices of a regular simplex inscribed in $S^{d}$ (strongly $1$-sharp), and two strongly $2$-sharp configurations (see \cite {DelGoeSei1977}): the Schl\"affi configuration of $N=27$ points on $S^5$ and the kissing (or McLaughlin) configuration of $N=275$ points on $S^{21}$. 
\begin {remark}\label {2.}
{\rm
Cardinality bounds for designs and $m$-distance sets by Delsarte, Goethals, and Seidel \cite {DelGoeSei1977} imply that strongly sharp configurations are tight designs. Combined with the result by Bannai and Dammerel \cite {BanDam1979}, these bounds imply immediately that there are no strongly $m$-sharp configurations on $S^d$ for $m\geq 3$ and $d\geq 2$.
}
\end {remark}

The following theorem holds.
\begin {theorem}\label {2m-1w}
Suppose $m,d\geq 1$, $f:[0,4]\to \RR$ is a function continuous on $[0,4]$, differentiable on $(0,4)$ with a concave derivative $f^{(2m-1)}$ on $(0,4)$. If $\omega_N=\{{\bf x}_1,\ldots,{\bf x}_N\}\subset S^{d}$ is a strongly $m$-sharp configuration, then the potential
$$
p_f({\bf x},\omega_N)=\sum\limits_{i=1}^{N}f(\left|{\bf x}-{\bf x}_i\right|^2),\ \ \ {\bf x}\in S^{d}, 
$$
attaints its absolute maximum over $S^{d}$ at every point of $\omega_N$.

If the concavity of $f^{(2m-1)}$ is strict on $(0,4)$, then $\omega_N$ contains all points of the absolute maximum of the potential $p_f(\cdot,\omega_N)$ over $S^{d}$.
\end {theorem}
A result similar to Theorem \ref {2m-1w} for absolute minima of strongly sharp configurations is established in~\cite {Borstiff}.

We next consider the class of antipodal (also called centrally symmetric) sharp configurations on the sphere.
A rather long list of examples of sharp configurations that occur in different dimensions is given in \cite {CohKum2007}. Among them, the antipodal ones are, for example, the configuration of $N=2m$ vertices of a regulat polygon inscribed in $S^1$ (which is $m$-sharp) and $N=2d+2$ vertices of a regular cross-polytope inscribed in $S^{d}$, which is $2$-sharp. 

Other examples are the configuration of $N=12$ vertices of a regular icosahedron inscribed in $S^2$, the kissing configuration of $N=56$ points on $S^6$, and the $552$-point configuration on $S^{22}$ (equiangular lines), which are all $3$-sharp. Also, the configuration of $N=240$ minimal non-zero vectors of the $E_8$ root lattice normalized to lie on the unit sphere $S^7$ and the kissing configuration of $N=4600$ vectors on $S^{22}$, which are both $4$-sharp. Finally, the configuration of $N=196560$ minimal non-zero vectors of the Leech lattice normalized to lie on $S^{23}$, which is $6$-sharp. We remark that any antipodal sharp configuration is a tight design (see \cite [Theorem 6.8]{DelGoeSei1977}).

We establish the following result.

\begin {theorem}\label {sharp}
Let $m,d\geq 1$, $f:[0,4]\to \RR$ be a function continuous on $[0,4]$, differentiable on $(0,4)$ with a convex derivative $f^{(2m-2)}$ on $(0,4)$. If $\omega_N=\{{\bf x}_1,\ldots,{\bf x}_N\}$ is an antipodal $m$-sharp configuration on the sphere $S^{d}$, then the potential
$$
p_f({\bf x},\omega_N)=\sum\limits_{i=1}^{N}f(\left|{\bf x}-{\bf x}_i\right|^2),\ \ \ {\bf x}\in S^{d}, 
$$
attaints its absolute maximum over $S^{d}$ at every point of $\omega_N$.

If the convexity of $f^{(2m-2)}$ is strict on $(0,4)$, then $\omega_N$ contains all points of absolute maximum of the potential $p_f(\cdot,\omega_N)$ over $S^{d}$.
\end {theorem}
The problem of finding the absolute minimum the potential of antipodal sharp configurations sometimes requires a more subtle approach. We deal with it in \cite {Borsymmetric,Borstiff}. 
We remark that in the case $m=1$ in Theorem \ref {sharp}, there is only one choice for the configuration $\omega_N$. Namely, $N=2$ and $\omega_N=\left\{{\bf x}_1,-{\bf x}_1\right\}$ for some ${\bf x}_1\in S^d$.

The proofs of Theorems \ref {2m-1w} and \ref {sharp} utilize the linear programming approach (see \cite {CohKum2007} and references therein) except here we do not need to show the positive definiteness of the interpolating polynomial.
We next state a solution to Problem \ref {mp} that follows from Theorems \ref {2m-1w} and \ref {sharp} and the universal optimality result from \cite {CohKum2007} for the minimal energy problem. 
\begin {theorem}\label {polarization}
Suppose $d\geq 1$ and $f:[0,4]\to \RR$ is a function continuous on $[0,4]$ and completely monotone on $(0,4]$ after possibly adding a constant. Suppose also that $\OL\omega_N$ is a strongly sharp or a sharp antipodal $N$-point configuration on $S^{d}$. Then 
\begin {equation}\label {O'}
Q_f(\OL \omega_N;S^{d})\leq Q_f(\omega_N;S^{d})
\end {equation}
for any $N$-point configuration (multiset) $\omega_N\subset S^{d}$. 

If, in addition, $f$ is strictly completely monotone on $(0,4]$ (after possibly adding a constant) and $\omega_N$ is such that equality holds in \eqref {O'}, then $\omega_N$ is sharp and the same dot products occur in $\omega_N$ and $\OL\omega_N$.
\end {theorem}
We remark that the case $N=d+2$ and $\OL\omega_N$ being the set of vertices of the regular $(d+1)$-simplex inscribed in $S^{d}$ in Theorem \ref {polarization} does not follow from the results of \cite {Borsimplex} on the maximal polarization for a simplex.

Theorem \ref {polarization} has the following two immediate consequences. The first one is for the Riesz potential with $-2<s<0$.
\begin {corollary}\label {R}
Suppose $d\geq 1$, $0<s<2$, and $\ \OL\omega_N=\{{\bf x}_1,\ldots,{\bf x}_N\}$ is a strongly sharp or a sharp antipodal configuration on $S^{d}$. Then
\begin {equation}\label {E1}
\min\limits_{{\bf x}\in S^{d}}\sum\limits_{i=1}^{N}\left|{\bf x}-{\bf y}_i\right|^s\leq \min\limits_{{\bf x}\in S^{d}}\sum\limits_{i=1}^{N}\left|{\bf x}-{\bf x}_i\right|^s
\end {equation}
holds for any configuration $\omega_N=\{{\bf y}_1,\ldots,{\bf y}_N\}\subset S^{d}$.

If $\omega_N$ is such that equality holds in \eqref {E1}, then $\omega_N$ is sharp and the same dot products occur in $\omega_N$ and $\OL \omega_N$.
\end {corollary}
We remark that for $N=3$, $d=1$, and $\OL\omega_N$ being the set of vertices of an equlateral triangle inscribed in $S^1$, Corollary \ref {R} was proved in \cite [p. 250]{Sto1975circle}.
The second consequence of Theorem \ref {polarization} is for Gaussian kernels.
\begin {corollary}\label {R'}
Suppose $d\geq 1$, $\sigma>0$, and $\OL\omega_N=\{{\bf x}_1,\ldots,{\bf x}_N\}$ is a strongly sharp or a sharp antipodal configuration on $S^{d}$. Then
\begin {equation}\label {E2}
\max\limits_{{\bf x}\in S^{d}}\sum\limits_{i=1}^{N}e^{-\sigma\left|{\bf x}-{\bf x}_i\right|^2}\leq \max\limits_{{\bf x}\in S^{d}}\sum\limits_{i=1}^{N}e^{-\sigma\left|{\bf x}-{\bf y}_i\right|^2}
\end {equation}
holds for any configuration $\omega_N=\{{\bf y}_1,\ldots,{\bf y}_N\}\subset S^{d}$.

If $\omega_N$ is such that equality holds in \eqref {E2}, then $\omega_N$ is sharp and the same dot products occur in $\omega_N$ and $\OL\omega_N$.
\end {corollary}
Theorem \ref {polarization} and Corollaries \ref {R} and \ref {R'} imply immediately the optimality of configurations mentioned in paragraphs preceding Remark \ref {2.} and Theorem~\ref {sharp} for the Riesz $s$-kernel ($-2<s<0$) and for the Gaussian kernel.

\section {Auxiliary results}

The following interpolation lemma is crucial for establishing Theorems \ref {2m-1w} and~\ref {sharp}. Let $\mathbb P_n$ denote the space of all polynomials of degree at most $n$.

\begin {lemma}\label {Interp}
Suppose $d,m\geq 1$, $-1\leq t_1<\ldots<t_m<t_{m+1}=1$ are arbitrary nodes, $\nu=0$ or $1$, and $L:=2m-2+\nu$. If $\nu=0$ we assume that $t_1=-1$ and, if $\nu=1$, we assume that $t_1>-1$. Suppose also that $g:[-1,1]\to \RR$ is a function continuous on $[-1,1]$ and differentiable on $(-1,1)$ such that the derivative $g^{(L)}$ is convex on $(-1,1)$. Let $q\in \mathbb P_{L+1}$ be the polynomial such that
\begin {equation}\label {i3}
\begin {split}
&q(t_i)=g(t_i),\ \ i=1,\ldots,m+1,\ \ q'(t_i)=g'(t_i),\ \ i=2,\ldots,m,\\
& \text{and, if}\ \  \nu=1,\ we\  additionally \ assume\  that \ \ q'(t_1)=g'(t_1).
\end {split}
\end {equation}
Then 
\begin {equation}\label {i5}
q(t)\geq g(t), \ \ t\in [-1,1]. 
\end {equation}
If, in addition, the derivative $g^{(L)}$ is strictly convex on $(-1,1)$, the inequality in \eqref {i5} is strict for $t\in [-1,1]\setminus \{t_1,\ldots,t_{m+1}\}$. 
\end {lemma}
\begin {proof}
The assertion of the lemma is trivial when $L=0$; i.e., when $m=1$ and $\nu=0$. Therefore, we assume that $L\geq 1$.
Inequality \eqref {i5} holds trivially for $t=t_1,\ldots,t_{m+1}$. We choose arbitrary $x\in [-1,1]\setminus \{t_1,\ldots,t_{m+1}\}$ and denote
$$
h(t):=q(t)+b(t-t_1)^{1+\nu}(t-t_{m+1})\prod_{i=2}^{m}(t-t_i)^2,
$$
where the constant $b$ is chosen so that $h(x)=g(x)$. Observe that the polynomial $h$ also satisfies interpolation conditions \eqref {i3}. The function $\tau(t):=g(t)-h(t)$ has at least $m+2$ distinct zeros in $[-1,1]$ at least $m-1+\nu$ of which have multiplicity two and are located in $(-1,1)$. By the Rolle's theorem, the derivative $\tau'$ has at least $2m+\nu=L+2$ distinct zeros in $(-1,1)$, which implies that the derivative $\tau^{(L)}$ has at least three distinct zeros in this interval. The polynomial $h^{(L)}$ has degree at most two with the coefficient of $t^2$ being $\frac {(L+2)!}{2}b$. It interpolates the convex function $g^{(L)}$ at at least three distinct points in $(-1,1)$. Then $b\geq 0$. Since $(x-t_1)^{1+\nu}>0$ both for $\nu=0$ and $\nu=1$, we have
\begin {equation}\label {i7}
g(x)=h(x)=q(x)+b(x-t_1)^{1+\nu}(x-t_{m+1})\prod_{i=2}^{m}(x-t_i)^2\leq q(x),
\end {equation}
which proves \eqref {i5}.
If $g^{(L)}$ is strictly convex on $(-1,1)$, then $b\neq 0$ (otherwise $g^{(L)}$ would be interpolated by the polynomial $h^{(L)}$ of degree at most one at three distinct points). Therefore, $b>0$ and the inequality in \eqref {i7} is strict for $x\in [-1,1]\setminus \{t_1,\ldots,t_{m+1}\}$. 
\end {proof}

Denote 
$$
w_{d}(t):=\gamma_{d}(1-t^2)^{d/2-1}, 
$$
where the constant $\gamma_{d}$ is chosen so that $w_{d}$ is a probability density on $[-1,1]$. We will need the following special case of the Funk-Hecke formula 
\begin {equation}\label{FnH}
\int_{S^{d}}g({\bf x}\cdot {\bf y})\ \! d\sigma_{d}({\bf x})=\int_{-1}^{1}g(t)w_{d}(t)\ \! dt,\ \ \ {\bf y}\in S^{d},
\end {equation}
where $g$ is any integrable function on $[-1,1]$ with the weight $w_{d}$.


We next recall the following property of spherical designs.
\begin {lemma}\label {const}
Let $d\geq 1$, $n\geq 0$, and $\omega_N=\{{\bf x}_1,\ldots,{\bf x}_N\}$ be a spherical $n$-design on $S^d$. Then for every polynomial $Q\in \mathbb P_n$, the sum 
$$
J_Q({\bf z}):=\sum\limits_{i=1}^{N}Q\({\bf z}\cdot{\bf x}_i \)
$$ 
is constant over ${\bf z}\in S^d$.
\end {lemma}
\begin {proof}
The assertion of Lemma \ref {const} is trivial when $n=0$.
Let $n\geq 1$ and ${\bf z}\in S^d$ be arbitrary. Since the polynomial $p({\bf x})=Q({\bf z}\cdot {\bf x})$ has degree at most $n$, using \eqref{FnH} we have
\begin {equation*}
\begin {split}
J_Q({\bf z})=\sum\limits_{i=1}^{N}Q\({\bf z}\cdot{\bf x}_i \)=N\int_{S^d}Q({\bf z}\cdot {\bf x})\ \! d\sigma_d({\bf x})=N\int_{-1}^{1}Q(t)w_d(t)\ \! dt,
\end {split}
\end {equation*}
which is independent of ${\bf z}$.
\end {proof}

%
%
%
To establish the uniqueness in Theorems \ref {2m-1w} and~\ref {sharp}, we need the following result.
\begin {lemma}\label {a5}
Let $d\geq 1$, $n\geq 0$, and $\omega_N\subset S^d$ be a $2n$-design. Then for every point ${\bf z}\in S^d$, the set $D({\bf z},\omega_N)$ contains at least $n+1$ distinct elements.
\end {lemma}
\begin {proof}
The assertion of the lemma is trivial for $n=0$. Let $n\geq 1$ and
assume to the contrary that there is a point ${\bf z}\in S^d$ with $D({\bf z},\omega_N)$ containing $k\leq n$ distinct elements. Denote them by $t_1<\ldots<t_k$ and let $p(t):=(t-t_1)^2\cdot \ldots \cdot (t-t_k)^2$. Since $p$ vanishes at each dot product ${\bf z}\cdot {\bf x}_i$ and $\omega_N$ is a $2n$-design, using the Funk-Hecke formula \eqref {FnH} we have
$$
0=\sum\limits_{i=1}^{N}p({\bf z}\cdot {\bf x}_i)=N\int_{S^d}p({\bf z}\cdot {\bf x})\ \! d\sigma_d({\bf x})=N\int_{-1}^{1}p(t)w_d(t)\ \! dt>0.
$$
This contradiction proves the assertion of the lemma.
\end {proof}

For the proof of Theorem \ref {2m-1w}, we also need the following statement.
\begin {lemma}\label {2m}
Let $d,m\in \NN$ and $\OL\omega_N$ be a strongly $m$-sharp configuration on $S^d$. Then no dot product between vectors in $\OL\omega_N$ equals $-1$. 
\end {lemma}
\begin {proof}
Assume to the contrary that one of the dot products between vectors in $\OL\omega_N$ equals $-1$. Let $-1=\tau_1<\tau_2<\ldots<\tau_m<1$ be all values of the dot product between distinct vectors in $\OL\omega_N$. Denote by $p(t)=(t+1)(t-\tau_2)^2\cdots(t-\tau_m)^2(t-1)$ and let ${\bf x}_1,\ldots,{\bf x}_N$ be points in $\OL\omega_N$. Since the polynomial $h({\bf x})=p({\bf x}_1\cdot {\bf x})$ has degree $2m$ and $\OL\omega_N$ is a $2m$-design, using the Funk-Hecke formula \eqref {FnH}, we have
$$
0=\sum\limits_{i=1}^{N}p({\bf x}_1\cdot {\bf x}_i)=N\int\limits_{S^d}p({\bf x}_1\cdot {\bf x})\ \!d\sigma_d({\bf x})=N\int\limits_{-1}^{1}p(t)w_d(t)\ \!dt<0.
$$
This contradiction shows that no dot product between vectors in $\OL\omega_N$ is $-1$.
\end {proof}

Finally, we establish a general statement used in
the proofs of Theorems \ref {2m-1w} and \ref {sharp}.
\begin {lemma}\label {gen}
Suppose that $d,m\geq 1$ and that a function $g:[-1,1]\to [-\infty,\infty)$ and a set $A$ consisting of numbers $-1\leq t_1<\ldots<t_{m+1}\leq 1$ are arbitrary. Suppose also that $q$ is a polynomial such that $q(t_i)=g(t_i)$, $i=1,\ldots,m+1$, and
\begin {equation}\label {1geq}
q(t)\geq g(t),\ \ t\in [-1,1].
\end {equation}
If $\omega_N=\{{\bf x}_1,\ldots,{\bf x}_N\}\subset S^{d}$ is a spherical $L$-design, where $L\geq {\rm deg}\ \! q$, and ${\bf y}\in S^{d}$ is any point such that $D({\bf y},\omega_N)\subset A$, then the potential 
$$
p^g({\bf x},\omega_N):=\sum\limits_{i=1}^{N}g({\bf x}\cdot {\bf x}_i), \ \ \ {\bf x}\in S^{d},
$$
attains its absolute maximum over $S^{d}$ at the point ${\bf y}$. 

If, in addition, the inequality in \eqref {1geq} is strict for every $t\in [-1,1]\setminus A$, then any point ${\bf z}\in S^{d}$ such that $D({\bf z},\omega_N)\not\subset A$ is not a point of absolute maximum of $p^g(\ \!\cdot\ \!,\omega_N)$ on~$S^{d}$.
\end {lemma}
\begin {proof}
By Lemma \ref {const}, since $\omega_N$ is a spherical $L$-design, the sum $\sum_{i=1}^{N}q({\bf x}\cdot {\bf x}_i)$ is constant over ${\bf x}\in S^d$.
Since $D({\bf y},\omega_N)\subset A$ and $q$ interpolates $g$ on $A$, taking into account inequality \eqref {1geq}, for arbitrary ${\bf x}\in S^d$, we have
\begin {equation}\label {1ineq}
p^g({\bf x},\omega_N)=\sum\limits_{i=1}^{N}g({\bf x}\cdot{\bf x}_i)\leq \sum\limits_{i=1}^{N}q({\bf x}\cdot {\bf x}_i)=\sum\limits_{i=1}^{N}q({\bf y}\cdot {\bf x}_i)=\sum\limits_{i=1}^{N}g({\bf y}\cdot {\bf x}_i)=p^g({\bf y},\omega_N),
\end {equation}
which shows that ${\bf y}$ is a point of absolute maximum.

Assume, in addition, that inequality \eqref {1geq} is strict for every $t\in [-1,1]\setminus A$. For any ${\bf z}\in S^d$ with $D({\bf z},\omega_N)\not\subset A$, there is an index $k$ such that ${\bf z}\cdot {\bf x}_k\notin A$ and, hence, $g({\bf z}\cdot {\bf x}_k)<q({\bf z}\cdot {\bf x}_k)$. Then the inequality in \eqref {1ineq} is strict for ${\bf x}={\bf z}$; that is ${\bf z}$ is not a point of absolute maximum of $p^g(\cdot,\omega_N)$.
\end {proof}

\section {Proof of Theorems \ref {2m-1w} and \ref {sharp}}

\begin {proof}[Proof of Theorem \ref {2m-1w}]
Let ${\bf y}\in \omega_N$ be arbitrary point and let $-1< \tau_1<\ldots<\tau_m<\tau_{m+1}=1$ be the distinct values of dot products between points in $\omega_N$ (we have $\tau_1>-1$ in view of Lemma \ref {2m}). By assumption, the function $g(t)=f(2-2t)$ has a convex derivative of order $2m-1$ on $(-1,1)$. Let $Q\in \mathbb P_{{2m}}$ be the polynomial such that $Q(\tau_i)=g(\tau_i)$, $i=1,\ldots, m+1$, and $Q'(\tau_i)=g'(\tau_i)$, $i=1,\ldots,m$. By Lemma \ref {Interp} (with $\nu=1$) we obtain that $Q(t)\geq g(t)$, $t\in [-1,1]$. Since $\omega_N$ is a $2m$-design and $D({\bf y},\omega_N)\subset \{\tau_1,\ldots,\tau_{m+1}\}$, by Lemma~\ref {gen} the potential
$$
p_f({\bf x},\omega_N)=\sum_{i=1}^{N}f(\left|{\bf x}-{\bf x}_i\right|^2)=\sum_{i=1}^{N}g({\bf x}\cdot {\bf x}_i)=p^g({\bf x},\omega_N)
$$
attains its absolute maximum over $S^d$ at ${\bf y}$.

If $f^{(2m-1)}$ is strictly concave then $g^{(2m-1)}$ is strictly convex. Lemma \ref {Interp} implies that $Q(t)>g(t)$ for $t\in [-1,1]\setminus \{\tau_1,\ldots,\tau_{m+1}\}$. Let ${\bf z}\in S^d\setminus\omega_N$ be arbitrary. By Lemma \ref {a5}, the set $D({\bf z},\omega_N)$ contains at least $m+1$ distinct values and does not contain $\tau_{m+1}=1$ (since ${\bf z}\notin\omega_N$). Then $D({\bf z},\omega_N)$ is not contained in $\{\tau_1,\ldots,\tau_{m+1}\}$. By Lemma \ref {gen}, the potential $p^g(\cdot,\omega_N)=p_f(\cdot,\omega_N)$ does not attain its absolute maximum at ${\bf z}$.
\end {proof}

\begin {proof}[Proof of Theorem \ref {sharp}]
Let ${\bf y}\in \omega_N$ be arbitrary point and let $-1=\tau_1<\ldots<\tau_m<\tau_{m+1}=1$ be distinct values of the dot product between vectors in $\omega_N$. By assumption, the function $g(t)=f(2-2t)$ has a convex derivative of order $2m-2$ on $(-1,1)$. Let $Q\in \mathbb P_{{2m-1}}$ be the polynomial such that $Q(\tau_i)=g(\tau_i)$, $i=1,\ldots, m+1$, and $Q'(\tau_i)=g'(\tau_i)$, $i=2,\ldots,m$. By Lemma \ref {Interp} (with $\nu=0$) we obtain that $Q(t)\geq g(t)$, $t\in [-1,1]$. Since $\omega_N$ is a $(2m-1)$-design, by Lemma~\ref {gen}, the potential
$
p_f({\bf x},\omega_N)=p^g({\bf x},\omega_N)
$
attains its absolute maximum over $S^d$ at ${\bf y}$.

If $f^{(2m-2)}$ is strictly convex then so is $g^{(2m-2)}$. Lemma \ref {Interp} implies that $Q(t)>g(t)$ for $t\in [-1,1]\setminus \{\tau_1,\ldots,\tau_{m+1}\}$. Let ${\bf z}\in S^d\setminus\omega_N$ be arbitrary. By Lemma \ref {a5} with $n=m-1$, the set $D({\bf z},\omega_N)$ contains at least $m$ distinct values. It does not contain $\tau_{m+1}=1$. It also does not contain $\tau_1=-1$ (if it did, then $-{\bf z}$ would be in $\omega_N$ and, since $\omega_N$ is antipodal, ${\bf z}$ itself would be in $\omega_N$). Thus, $D({\bf z},\omega_N)\not\subset\{\tau_1,\ldots,\tau_{m+1}\}$. By Lemma \ref {gen}, the potential $p_f(\cdot,\omega_N)=p^g(\cdot,\omega_N)$ does not attain its absolute maximum at ${\bf z}$.
\end {proof}

\section {Proof of Theorem \ref {polarization}}

We use the following well-known result from \cite {CohKum2007}.
\begin {theorem}\label {CohKum}
Let $f : (0, 4] \to \RR$ be completely monotone, and let $\mathcal C \subset S^{d}$ be a
sharp configuration. If $\mathcal C' \subset S^{d}$ is any subset
satisfying $\#\mathcal C'
 = \#\mathcal C$, then
\begin {equation}\label {1.1}
\sum\limits_{
{\bf x},{\bf y}\ \!\!\in \ \!\!\mathcal C' ,\ \! {\bf x}\neq {\bf y}}
f\(|{\bf x} -{\bf  y}|^2\)\geq  \sum\limits_{
{\bf x},{\bf y}\ \!\!\in \ \!\!\mathcal C,\ \! {\bf x}\neq {\bf y}}
f\(|{\bf x} - {\bf y}|^2\).
\end {equation}
If $f$ is strictly completely monotone, then equality in \eqref {1.1} implies that $\mathcal C'$ is also a sharp configuration and the same distances occur in $\mathcal C$ and $\mathcal C'$.
\end {theorem}

\begin {proof}[Proof of Theorem \ref {polarization}]
Without a loss of generality, we can assume that $f$ is completely monotone on $(0,4]$ (otherwise we add an appropriate contstant to~$f$). The derivative $f^{(2m-1)}$ is concave on $(0,4)$ and the derivative $f^{(2m-2)}$ is convex on $(0,4)$. If $\OL\omega_N=\{{\bf x}_1,\ldots,{\bf x}_N\}$ is strongly $m$-sharp for some $m\geq 1$, we apply Theorem \ref {2m-1w}. If $\OL\omega_N$ is an $m$-sharp antipodal configuration for some $m\geq 1$, we apply Theorem~\ref {sharp}. In both cases, we conclude that the potential $p_f(\cdot,\OL\omega_N)$ achieves its absolute maximum over $S^{d}$ at points of $\OL\omega_N$. Let $\omega_N=\{{\bf y}_1,\ldots,{\bf y}_N\}$ be an arbitrary configuration (multiset) on $S^d$. By Theorem \ref {CohKum}, 
\begin {equation}\label {CohnK}
\sum\limits_{i=1}^{N}\sum\limits_{j=1\atop j\neq i}^{N}f\(\left|{\bf x}_i-{\bf x}_j\right|^2\)\leq \sum\limits_{i=1}^{N}\sum\limits_{j=1\atop j\neq i}^{N}f\(\left|{\bf y}_i-{\bf y}_j\right|^2\)
\end {equation}
whenever ${\bf y}_i\neq {\bf y}_j$ for $i\neq j$ in $\omega_N$.
Inequality \eqref {CohnK} still holds if ${\bf y}_i={\bf y}_j$ for some $i\neq j$, since in this case there is a sequence of $N$-point configurations on $S^d$ where points with different indices are distinct that converges to $\omega_N$.
For any fixed index $\ell$, we now have
\begin {equation*}
\begin {split}
Q_f(\OL\omega_N;S^{d})&=\max\limits_{{\bf x}\in S^{d}}\sum\limits_{i=1}^{N}f\(\left|{\bf x}-{\bf x}_i\right|^2\)=\sum\limits_{i=1}^{N}f\(\left|{\bf x}_\ell-{\bf x}_i\right|^2\)\\
&=\frac {1}{N}\sum\limits_{j=1}^{N}\sum\limits_{i=1}^{N}f\(\left|{\bf x}_j-{\bf x}_i\right|^2\)\\
&=\frac {1}{N}\(\sum\limits_{j=1}^{N}\sum\limits_{i=1\atop i\neq j}^{N}f\(\left|{\bf x}_j-{\bf x}_i\right|^2\)+N f(0)\).
\end {split}
\end {equation*}
From \eqref {CohnK} we have
\begin {equation}\label {Q}
\begin {split}
Q_f(\OL\omega_N;S^{d})&\leq \frac {1}{N}\(\sum\limits_{j=1}^{N}\sum\limits_{i=1\atop i\neq j}^{N}f\(\left|{\bf y}_j-{\bf y}_i\right|^2\)+N f(0)\)\\
&=\frac {1}{N}\sum\limits_{j=1}^{N}\sum\limits_{i=1}^{N}f\(\left|{\bf y}_j-{\bf y}_i\right|^2\)=\frac {1}{N}\sum\limits_{j=1}^{N}p_f({\bf y}_j,\omega_N)\\
&\leq \max\limits_{{\bf x}\in S^{d}}p_f({\bf x},\omega_N)=Q_f(\omega_N;S^{d}),
\end {split}
\end {equation}
which completes the proof of \eqref {O'}. 

Assume now that $f$ is strictly completely monotone on $(0,4]$ (after possibly adding a constant). Then equality in \eqref {O'} implies equality between the energies of $\OL\omega_N$ and $\omega_N$ in \eqref {CohnK}. If ${\bf y}_i\neq {\bf y}_j$ for all $i\neq j$, Theorem \ref {CohKum} immediately implies that $\omega_N$ is sharp and the same distances occur in $\omega_N$ and $\OL\omega_N$. 

It remains to show that when ${\bf y}_i={\bf y}_j$ for some $i\neq j$, the configuration $\omega_N$ is not energy minimizing. 
Denote by $-1\leq t_1<\ldots<t_m<1$ the dot products that occur between distinct points in $\OL\omega_N$. The function $g(t):=f(2-2t)$ is continuous on $[-1,1]$ and strictly absolutely monotone on $[-1,1)$; that is, $g^{(k)}>0$ on $(-1,1)$ for every $k\geq 0$.
Let $p\in \mathbb P_{2m-1}$ be the Hermite interpolating polynomial for $g$ at $t_1,\ldots,t_m$. Then (cf., e.g., \cite [Theorem 3.5.1]{Dav1963})
\begin {equation}\label {D}
g(t)\geq p(t), \ \ t\in [-1,1],
\end {equation}
and the inequality in \eqref {D} is strict for $t\neq t_1,\ldots,t_m$. Since one of the dot products occuring in $\omega_N$ equals $1$ and $g(1)>p(1)$, we have 
\begin {equation}\label {D1}
\sum_{1\leq i\neq j\leq N}f(\left|{\bf y}_i-{\bf y}_j\right|^2)=\sum_{1\leq i\neq j\leq N}g({\bf y}_i\cdot{\bf y}_j)>\sum_{1\leq i\neq j\leq N}p({\bf y}_i\cdot {\bf y}_j).
\end {equation}
Then (see the proof in \cite [Sections $2-6$]{CohKum2007})
$$
\sum_{1\leq i\neq j\leq N}p({\bf y}_i\cdot{\bf y}_j)\geq \sum\limits_{1\leq i\neq j\leq N}{p({\bf x}_i\cdot{\bf x}_j)}=\sum\limits_{1\leq i\neq j\leq N}{g({\bf x}_i\cdot{\bf x}_j)}=\sum\limits_{1\leq i\neq j\leq N}{f(\left|{\bf x}_i-{\bf x}_j\right|^2)},
$$
which combined with \eqref {D1} shows that $\omega_N$ is not energy minimizing. Then the first inequality in \eqref {Q} is strict and, hence, $\omega_N$ does not satisfy the equality in~\eqref {O'}.
\end {proof}

\end {document}